\documentclass[]{article}  
\pdfoutput=1
\usepackage{graphicx}
\usepackage{float}
\usepackage{amsthm}
\usepackage{amssymb}
\usepackage{amsmath}
\usepackage{mathrsfs}
\usepackage{indentfirst}
\usepackage{geometry}
\usepackage{authblk}
\usepackage{blkarray}
\renewcommand\figurename{Fig}
\newtheorem{Theorem}{Theorem}[section]
\theoremstyle{definition}

\newtheorem{Lemma}[Theorem]{Lemma}

\theoremstyle{definition}

\title{A characterization of graphs $G$ with nullity $n(G)-d(G)-1$}
\author{Songnian Xu\thanks{Corresponding author. E-mail address: xsn131819@163.com.}}
\affil{\textit{School of Mathematics, China University of Mining and Tecnology, Xuzhou, China.}}
\date{}

\begin{document}
\baselineskip 17pt

\title{A characterization of graphs $G$ with nullity $n(G)-d(G)-1$}

\author{Songnian Xu\\
{\small  Department of Mathematics, China University of Mining and Technology}\\
{\small Xuzhou, 221116, P.R. China}\\
{\small E-mail: xsn1318191@cumt.edu.cn}\\ \\
Wang\thanks{Corresponding author}\\
{\small Department of Mathematics, China University of Mining and Technology}\\
{\small Xuzhou 221116, P.R. China}\\
{\small E-mail: }}

\date{}
\maketitle

\begin{abstract}
For a connected graph $G$ with order $n$, let $e(G)$ represent the number of its distinct eigenvalues, and let $d$ denote its diameter.
We denote the eigenvalue multiplicity of $\mu$ in $G$ by $m_G(\mu)$.
It is well established that the inequality $e(G) \geq d + 1$ implies that when $\mu$ is an eigenvalue of $P_{d+1}$, it follows that $m_G(\mu) \leq n - d$; otherwise, for any real number $\mu$, we have $m_G(\mu) \leq n - d - 1$.
A graph is termed minimal if $e(G) = d + 1$.
In 2013, Wong et al. characterized all minimal graphs for which $m_G(0) = n - d$.
In this article, we provide a complete characterization of the graphs $G$ such that $m_G(0) = n - d - 1$.
\end{abstract}

\let\thefootnoteorig\thefootnote
\renewcommand{\thefootnote}{\empty}

\section{Introduction}
In this paper, we consider only simple, connected and finite graphs.
A simple undireted graph $G$ is denoted by $G=(V(G),E(G))$, where $V(G)$ is the vertex set and $E(G)$ is the edge set.
A graph $H$ is called a subgraph of a graph $G$ if $V(H)\subseteq V(G)$ and $E(H)\subseteq E(G)$.
Further, if any two vertices of $V(H)$ are adjacent in $H$ if and only if they are adjacent in $G$, we say that $H$ is an induced subgraph of $G$ and denote this relation as $H\leq G$.
The order of $G$ is the number of vertices of $G$, denoted by $|G|$.
For $x\in V(G)$ and $H\leq G$, $N_H(x)=\{u\in V(H) |\ uv\in E(G)\}$.
Let $K\subseteq V(G)$, traditionally, the subgraph of $G$ induced by $K$, written as $G[K]$.
And we sometimes write $G-K$ or $G-G[K]$ to denote $G[V(G)\backslash V(K)]$.
A vertex $v$ of $G$ is said to be pendant if $d(v)=1$, where $d(v)$ denotes the number of adjacent vertices of $v$ in $V(G)$.
For $x, y \in V(G)$, distance $d(x, y)$ represents the length of the shortest path between $x$ and $y$, and $d(x, H)=min\{d(x,y)|y\in V(H)\}$ for $H\leq G$.
We denote by $K_{m,n}$ the complete bipartite graph.
Let $A(G)$ denote the adjacency matrix of graph $G$, which is a square matrix and $a_{ij}=1$ if and only if $v_i\sim v_j$, otherwise $a_{ij}=0$.

Let $u$, $v\in V(G)$ and $N_G(u) = N_G(v)$, then we call $u$ and $v$ are twins in $G$.
The graph $G$ is called reduced if $G$ contains no twin vertices.
We denote the reduced graph corresponding to $G$ by $G^r$.

Let $m_G(\mu)$ denote the multiplicity of the eigenvalue $\mu$ in the adjacency matrix $A(G)$.
In particular, the multiplicity of $0$ is denoted by $\eta(G)=m_G(0)$.
In 1957, Collatz and Sinogowitz \cite{CL1} first raised the problem of characterizing all singular graphs.
After then, nullities of various graphs have attracted much attention of researchers because of the relationship between the nullity of a graph representing an alternate hydrocarbon (a bipartite molecular graph) and its molecular stability \cite{LHC1}.
Moreover, several papers have concentrated on exploring the relationships between the nullity of graphs and certain structural parameters, such as pendant vertices and matching numbers, among others (see \cite{CQQ1}, \cite{GJM1}, \cite{LX1}, \cite{MX1}, \cite{MX2}, \cite{RS1}, \cite{SY1}, \cite{WL1}, \cite{WL2}).
Furthermore, the investigation of nullity in various types of graphs remains a prominent topic in graph theory, including oriented graphs, signed graphs, and mixed graphs (see \cite{CC1}, \cite{CQQ2}, \cite{HS1}, \cite{HS2}, \cite{Wong2}, \cite{XF1}).

In 2013, Wong \cite{Wong1} characterize graphs with maximum diameter among all connected graph with rank $n$.
In 2022, Wang \cite{Wang1} characterized the graphs for which $m_G(0) = n - d$ using the method of star complements.
In the same year, Chang and Li \cite{Chang1} characterized the graphs for which $m_G(0) = n - g - 1$.
In 2023, Du and Fonseca \cite{Du1} extend the definition of minimal graphs on
adjacency matrices to real symmetric matrices. They characterize all the trees for which
there is a real symmetric matrix with nullity $n-d$ and $n-d-1$.
Therefor, we investigate the situation $m_G(0)=n-d-1$.

\section{Preliminaries}
\begin{Lemma}\cite{AB}
Let $v$ be a vertex of $G$, then $m_G(\mu)-1\leq m_G(\mu)\leq m_G(\mu)+1$
\end{Lemma}

The eigenvalues of $P_n$ are $\{2cos\frac{i\pi}{n+1}|i = 1, 2, ..., n\}$ and the eigenvalues of $C_n$ are $\{2cos\frac{2i\pi}{n}|i =0, 1, 2, \ldots, n-1\}$, we have the following lemma.

\begin{Lemma}
$m_{P_{n+1}}(0)\leq1$ and $m_{C_n}(0)\leq 2$.
Further, $m_{P_{n+1}}(0) = 1$ if and only if $n \equiv 0 \ (\text{mod} \ 2)$, and $m_{C_n}(0) = 2$ if and only if $n \equiv 0 \ (\text{mod} \ 4)$.
\end{Lemma}

\begin{Lemma}
If  $m_G(0) = n - d - 1$ and $P_{d+1}$ is the diameter path of $G$ with $P_{d+1} \leq H \leq G$, then it follows that $rank(A(H)) \geq rank(A(G)) - 1$.
\end{Lemma}
\begin{proof}
If 0 is not an eigenvalue of $P_{d+1}$, then by the interlacing theorem and $\eta(G)=n(G)-rank(G)$, we have $rank(H)= rank(P_{d+1})$.
Similarly, if 0 is an eigenvalue of $P_{d+1}$, then $rank(H) = rank(P_{d+1})$ or $rank(P_{d+1})+1$.
This completes the proof.
\end{proof}

\begin{Lemma}
Let $H$ be an induced subgraph of a connected graph $G$ for which $rank(A(H))\geq rank(A(G))-1$, and $v\in V(G)\setminus V(H)$.
If $v\nsim h$ and $N_H(v) = N_H(h)$ for a vertex $h\in V (H)$, then $N_G(v) = N_G(h)$.
\end{Lemma}
\begin{proof}
Assume there exists $x \in V(G)\setminus V(H)$ such that $x \in N_G(h)$ but $x \notin N_G(v)$.
We consider the adjacency matrix $A_1$ corresponding to the subgraph $H + x + v$.

$A_1=\bordermatrix{
& v & x& h \cr
 & 0& 0 & 0 &\alpha^{T} \cr
 &0&  0& 1 &\beta^{T} \cr
  &0& 1 & 0 &\alpha^{T} \cr
 &\alpha &\beta  &\alpha  &A(H-h)
}
$,

where $\alpha$, $\beta$ are $|H|-1$ dimensional vectors.
Let
\\

$Q_1=
\begin{pmatrix}
  1& -1 & 0 &0 \\
 0&  1& 0& 0 \\
  0& 0 & 1 &0 \\
 0 &0  &0  &I
\end{pmatrix}
$,
$Q_2=
\begin{pmatrix}
  1& 0 & 0 &0 \cr
 0&  1& 0 &0 \cr
  1& 0 & 1 &0 \cr
 \beta &0  &0  &I
\end{pmatrix}
$.
\\

Then

$Q_2Q_1A_1{Q_1}^{T}{Q_1}^{T}=\bordermatrix{
& v & x& h \cr
 & 0& -1& 0 &0 \cr
 &-1&  0& 0 &0 \cr
  &0& 0 & 0 &\alpha^{T} \cr
 &0 &0 &\alpha  &A(H-h)
}=$
$\bordermatrix{
& v & x&  \cr
 & 0& -1& 0  \cr
 &-1&  0& 0  \cr
 &0 &0  &A(H)
}=B_1$.
\\

Therefore, we have $\text{rank}(A_1) = \text{rank}(B_1) = \text{rank}(A(H)) + 2 \geq \text{rank}(A(G)) + 1$, which is evidently a contradiction.
Thus, it follows that $N_G(h) \subseteq N_G(v)$.
Similarly, by symmetry, let $H' = H - h + v$, which leads to $N_G(v) \subseteq N_G(h)$.
Therefore, we conclude that $N_G(v) = N_G(h)$.
\end{proof}

\begin{Lemma}
Let $H$ be an induced subgraph of a connected graph $G$ for which $rank(A(H))\geq rank(A(G))-1$, and $u, v\in V (G)\backslash V (H)$.
If $u\nsim v$ and $N_H(u) = N_H(v)$, then $N_G(u)= N_G(v)$.

\end{Lemma}
\begin{proof}
Let $H' = H + u$.
We have $rank(A(H')) = rank(A(H + u)) \geq rank(A(H)) \geq rank(A(G)) - 1$.
Since $N_{H'}(u) = N_{H'}(v)$, by Lemma 2.4, we can deduce that $N_G(u) = N_G(v)$.

\end{proof}

\begin{Lemma}\cite{FYZ1}
Let $G$ be a graph.
If two vertices $u$, $v$ satisfy $N_G(u)=N_G(v)$, then $m_G(0)=m_{G-u}(0)+1=m_{G-v}(0)+1$.
\end{Lemma}

\begin{Lemma}\cite{CDM1}
Let $G$ be a graph obtained from a given graph $H$ and $P_2$, disjoint with $H$, by identifying one pendant vertex of $P_2$ with a vertex $w$ of $H$. Then $m_G(0)=m_{G-w}(0)$
\end{Lemma}

\begin{Lemma}
Let $G^{r}$ be the reduced graph of a graph $G$.
Then $m_G(0)=n(G)-d(G)-1$ if and only if $m_{G^{r}}(0)=n(G^{r})-d(G^{r})-1$.
\end{Lemma}

\begin{proof}
First, it can be easily observed that $d(G^{r})=d(G)$.
Meanwhile, by Lemma 2.6, $m_G(0)=m_{G^{r}}(0)+n(G)-n(G^{r})$.
So, $m_G(0)=n(G)-d(G)-1$ if and only if $m_{G^{r}}(0)+(n(G)-n(G^{r})=n(G)-d(G)-1=n(G)-d(G^{r})-1$.
Then we have $m_G(0)=n(G)-d(G)-1$ if and only if $m_{G^{r}}(0)=n(G^{r})-d(G^{r})-1$.
\end{proof}

By Lemmas 2.8, we know that for the eigenvalue $0$, $G$ and $G^{r}$ possess identical properties.
Therefore, in this article, we consider only the case where the graph $G$ is reduced.

\section{Characterizations of graphs with $m_G(0)=n-d-1$ }

\begin{figure}[H]
  \centering
  \includegraphics[width=1.0\linewidth]{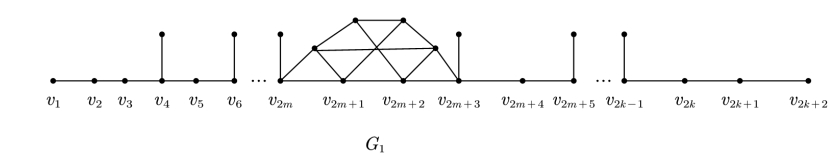}

\end{figure}

When  $m_G(0) = n - d - 1$ and $0$ is not an eigenvalue of $P_{d+1}$, we have $\text{rank}(G) = d + 1$.
This extremal graph has been completely characterized in \cite{WL3}.
Therefore, in the following article, we will only consider the case where $m_{P_{d+1}}(0) = 1$, which corresponds to $d \equiv 0 \ (\text{mod} \ 2)$.

\begin{Theorem}\cite{WL3}
Let $G$ be a connected graph with $d(G) = 2k + 1$.
Then, we have $d(G) + 1 \leq rank(G)$,
and equality holds if and only if $P_{d+1}\leq G \leq G_1$.
\end{Theorem}

First, we will present some graphs that will be needed for this study.

\begin{figure}[H]
  \centering
  \includegraphics[width=1.0\linewidth]{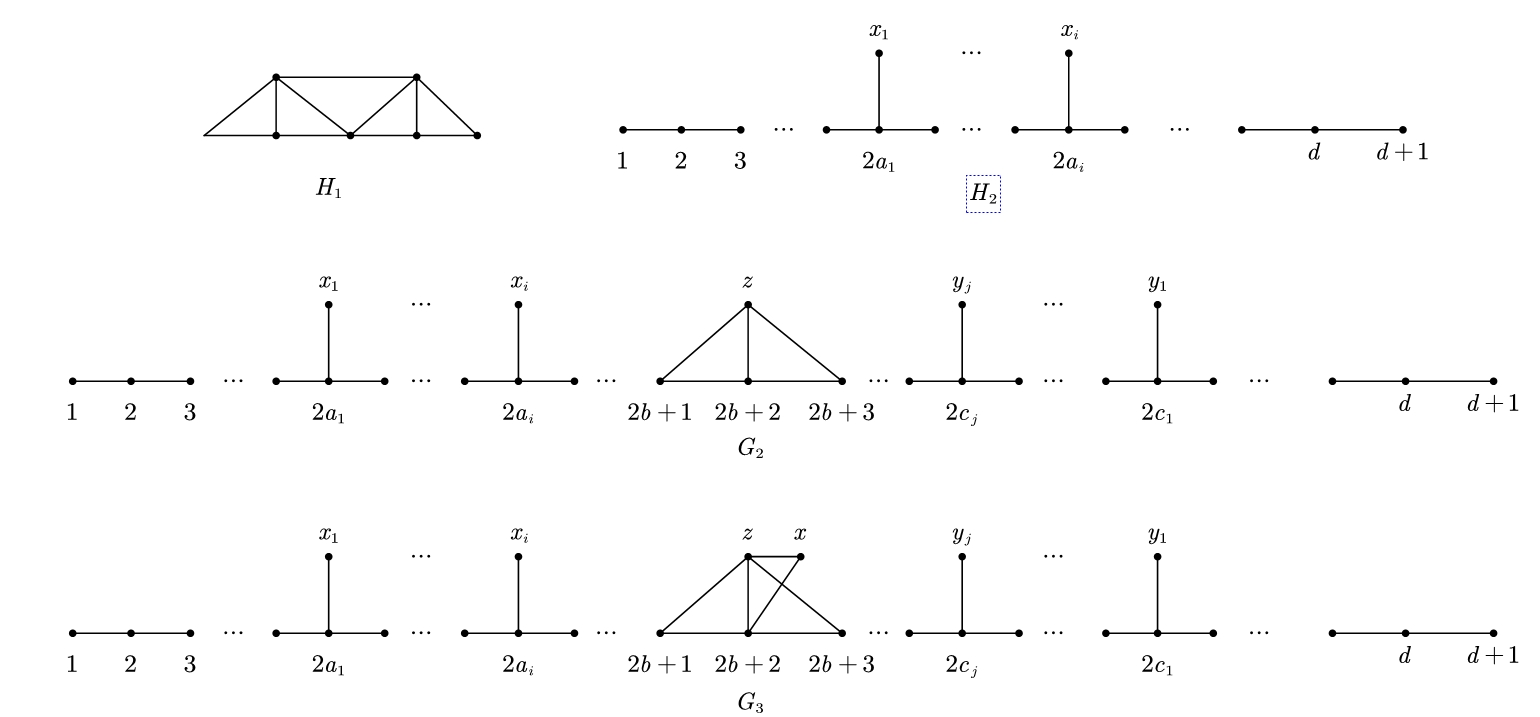}

\end{figure}

\begin{Theorem}
Let $G$ be a connected graph and $P = P_{d+1}=v_1\sim v_2\sim\cdots \sim v_{d+1}$ be the diameter path of $G$.
If $m_P(0)=1$ and $m_G(0)=n-d-1$, then $G\cong G_2$ or $G_3$.
\end{Theorem}

\begin{proof}

\textbf{Claim 1}: If $x\in V(G)\setminus V(P)$, then $d(x,P)=1$ .

If there exists an $x\in V(G)$ such that $d(x, P) = 2$, without loss of generality, let us assume $x \sim y \sim P$.
According to Lemma 2.7, we know that $m_{P+x+y}(0) = m_P(0) = 1$.
From Lemma 2.1, it follows that $m_G(0) \leq 1 + n - d - 3 = n - d - 2$, which leads to a contradiction.

Now, we define $V(G) \setminus V(P) = \bigcup_{i\geq 1} N_i$, where $N_i = \{ x : |N_{P}(x)| = i \}$.

\textbf{Claim 2}:  $d \equiv 0 \ (\text{mod} \ 2)$; if $x \in N_1$ and $x \sim v_m$, then $m \equiv 0 \ (\text{mod} \ 2)$.

First, we note that $1 \leq m_{P+x}(0)$; otherwise, by Lemma 2.1, we would obtain $m_G(0) < n - d - 1$.
Additionally, we also know that $m_{P+x}(0) \leq 2$ by $m_P(0)=1$.
For the general case, suppose $x_1, x_2, \ldots, x_m \in N_1 $.
We have $m \leq m_{P+x_1+x_2+\ldots+x_m} \leq m + 1$.

Since 0 is an eigenvalue of $P$, by Lemma 2.2, we conclude that $d \equiv 0 \ (\text{mod} \ 3)$.
If $m \equiv 1 \ (\text{mod} \ 2)$, then by applying Lemma 2.7, we find that $m_{P+x}(0) = m_{P_2}(0) = 0$, leading to a contradiction.
If $m \equiv 0 \ (\text{mod} \ 2)$, then $m_{P+x}(0) = 2$, which satisfies the condition.

\textbf{Claim 3}: If $x, y\in N_1$, then $x\nsim y$.

By Claim 2, we can assume without loss of generality that $x \sim v_{2a}$ and $y \sim v_{2b}$ with $a \leq b$.

If $a = b$ and $x \nsim y$, then $N_P(x) = N_P(y)$.
According to Lemma 2.4, this implies $N_G(x) = N_G(y)$, leading to a contradiction.
If $a = b$ and $x \sim y$, then by applying Lemma 2.7, we obtain $m_{P+x+y}(0) = m_{P_2}(0) + 1 = 1 < 2$, which is also a contradiction.

If $b > a$ and $x \nsim y$, then $m_{P+x+y}(0) = 3$, satisfying the required condition.
If $b > a$ and $x \sim y$, it follows that $b = a + 1$; otherwise, $v_1 v_2 \ldots v_{2a} xy v_{2b} \ldots v_{d+1}$ would form a path shorter than $P$, contradicting the assertion that $P$ is a diameter path.
When $b = a + 1$, we have $m_{P+x+y}(0) = 1 < 2$, which again leads to a contradiction.

\textbf{Claim 4}: $N_2=\emptyset$.

If there exists $u \in N_2$ such that $u \sim v_m$ and $u \sim v_t$ with ( $t > m$ ), and considering that $P$ is a diameter path, it follows that $t = m + 1$ or $t = m + 2$; otherwise, $v_1 v_2 \ldots v_t u v_m \ldots v_{d+1}$ would constitute a path shorter than $P$.

If $t = m + 2$, then $N_P(v_{m+1}) = N_P(u)$.
By Lemma 2.4, this implies $N_G(v_{m+1}) = N_G(u)$, leading to a contradiction.

If $t = m + 1$ and $m \equiv 0\ (\text{mod} \ 2)$ (or $m \equiv 1\ (\text{mod} \ 2)$), applying Lemma 2.7 gives us $m_{P+u}(0) = m_{P_2}(0) = 0$, which is also a contradiction.
Therefore, we conclude that $N_2=\emptyset$.

\textbf{Claim 5}: If $z\in N_3$, then $z \sim v_m, v_{m+1}, v_{m+2}$ and $m \equiv 1\ (\text{mod} \ 2)$.

For $z \in N_3$, suppose $z \sim v_{i_1}, v_{i_2}, v_{i_3}$, where we can assume $i_1 < i_2 < i_3$.
Considering that  $P$ is a diameter path, we have $i_3 = i_2 + 1 = i_1 + 2$.

If $i_1 \equiv 1\ (\text{mod} \ 2)$, then by Lemma 2.7, we have $m_{P+z}(0) = m_{P_2}(0) = 0$, which leads to a contradiction.

If $i_1 \equiv 1\ (\text{mod} \ 2)$, then $m_{P+z}(0) = 1 + m_{C_3}(0) = 1$, satisfying the required condition.

\textbf{Claim 6}: $|N_3|=1$.

If $|N_3| \geq 2 $, let $z, w \in N_3$.
According to Claim 5, we assume $z \sim v_{2a+1}, v_{2a+2}, v_{2a+3}$ and $w \sim v_{2b+1}, v_{2b+2}, v_{2b+3}$ where $b \geq a$.

If $b = a$ and $z \nsim w$, then $N_P(z) = N_P(w)$.
By Lemma 2.4, this implies $N_G(z) = N_G(w)$, leading to a contradiction.

If $b = a$ and $z \sim w$, then by Lemma 2.7, we have $m_{P + z + w}(0) = 1 + m_{K_4}(0) = 1$, which is also a contradiction.

If $b > a$ and $z \nsim w$, then by Lemma 2.1, we have $m_{P + z + w}(0) \leq 1 + m_{P + z + w - v_{2a+3}}(0) = 1$, leading to a contradiction.

If $b > a$ and $z \sim w$, considering that $P$ is a diameter path, it follows that $b = a + 1$.
In this case, $m_{P + z + w}(0) = m_{H_1}(0)$.
By computation, we can obtain $m_{H_1}(0) = 1$, which leads to $m_{P + z + w}(0) = 1<2$, which is also a contradiction.
Thus, $|N_3| \leq 1$.
If $|N_3| = \emptyset$, based on previous discussions, we understand that $G\cong H_2$.
According to Lemma 2.7, we find that $m_{G_3} = i + 1 = n - d > n - d - 1$, which presents another contradiction.
In summary, we conclude that $|N_3| = 1$.

\textbf{Claim 7}: Let $x \in N_1$ and $z \in N_2$, with $x \sim v_{2a}$ and $z \sim v_{2b+1}, v_{2b+2}, v_{2b+3}$.
Then, $ x \sim z$ if and only if $a = b + 1 $.

By symmetry, we may assume $a \leq b + 1$.
If $x \sim z$, considering that $P$ is a diameter path, it suffices to examine the cases where $a = b$ and $a = b + 1$.

If $a = b$, then $m_{P + x + z}(0) = 1 < 2$, leading to a contradiction.

If $a = b + 1$, then $m_{P + x + z}(0) = 2$, which satisfies the condition.
Conversely, if $a = b + 1$ and $x \nsim z$, then $m_{P + x + z}(0) = 1 < 2$, resulting in another contradiction.

In summary, $G_2$ and $G_3$ are graphs that satisfy all the conditions.
Furthermore, by applying Lemmas 2.6 and 2.7, we can determine that  $m_{G_2} = i + j + 1 = n - d - 1$ and $m_{G_3} = i + j + 2 = n - d - 1$. Thus, $G\cong G_2$ or $G_3$.
\end{proof}

\end{document}